\documentclass[11pt]{article}

\newtheorem{theorem}{Theorem}

 \RequirePackage{ifthen}
\RequirePackage{times}
\RequirePackage{mathptm}
\RequirePackage{pifont}

\usepackage{amssymb}
\usepackage{amsmath}
\usepackage{graphicx}

\usepackage{epsfig}

\usepackage{fullpage}

\newcommand\Li{\mathrm{Li}}

\begin{document}

\title{
Partition Polynomials: Asymptotics and Zeros
}

\author{Robert  P. Boyer, William M. Y. Goh}

\maketitle

\begin{abstract}
Let $F_n(x)$ be the partition polynomial
$\sum_{k=1}^n p_k(n) x^k$ where
 $p_k(n)$ is the number of partitions of $n$ with $k$ parts.
 We emphasize the computational experiments using
 degrees up to $70,000$ to
 discover the asymptotics of these polynomials.
 Surprisingly, the asymptotics of $F_n(x)$ 
 have two scales of orders $n$ and $\sqrt{n}$
and in three different regimes inside the unit disk.
 Consequently, the zeros converge to network of curves
 inside the unit disk given in terms of the dilogarithm.
   \end{abstract}

 \section{Introduction}\label{section:intro}
 
 Let $p(n)$ denote the number of partitions of a positive integer $n$;
 that is, the number of ways of writing $n$ additively. Euler identified their generating function:
 \[
 P(u)= \prod_{k=1}^\infty \frac{1}{1-u^k} = \sum_{k=1}^\infty p(n)u^k.
 \]
By inverting this equation, we find an integral form for $p(n)$ as
\[
p(n)= \frac{1}{2 \pi i} \oint_{|u|=r} \frac{ P(u)}{u^{n+1}} \, du, \quad
0<r<1.
\]
In their celebrated work, Hardy and Ramanujan in 1917 
\cite{hardy}
discovered the asymptotics of $p(n)$ by developing a new approach ``the circle method" to handle the asymptotics of such contour integrals which have dense singularities over the unit circumference. A basic form of their asymptotics is
\[
p(n) = \frac{1}{4n \sqrt{3}} e^{ \pi \sqrt{2n/3}}\,
\left(
1 + O \left( 1 / \sqrt{n} \right)
\right)
\]
We emphasize that Hardy and Ramanujan were guided by an extensive list 
that appears in their paper
of values of $p(n)$ up to $n=200$.

The {\sl partition polynomial}
 $F_n(x)$ is a refinement of the partition numbers:
\[
F_n(x) = \sum_{k=1}^n p_k(n) x^k
\]
where $p_k(n)$ is the number of partitions of $n$ with exactly $k$ parts. For example, $F_5(x)=x^5+x^4+2x^3+2x^2+x$. 
We note that the partition numbers are recovered as 
$F_n(1)=p(n)$. Further, $p_k(n)$ can be computed via the recurrence
$p_k(n)  =  p_{k-1}(n-1)  + p_k(n-k) $.

Richard Stanley plotted the zeros of $F_{200}$ and asked what happens as their degrees go to infinity. See Figure 1(a).

   \begin{figure}\begin{center}
\includegraphics[height=6.5cm,width=6.5cm]{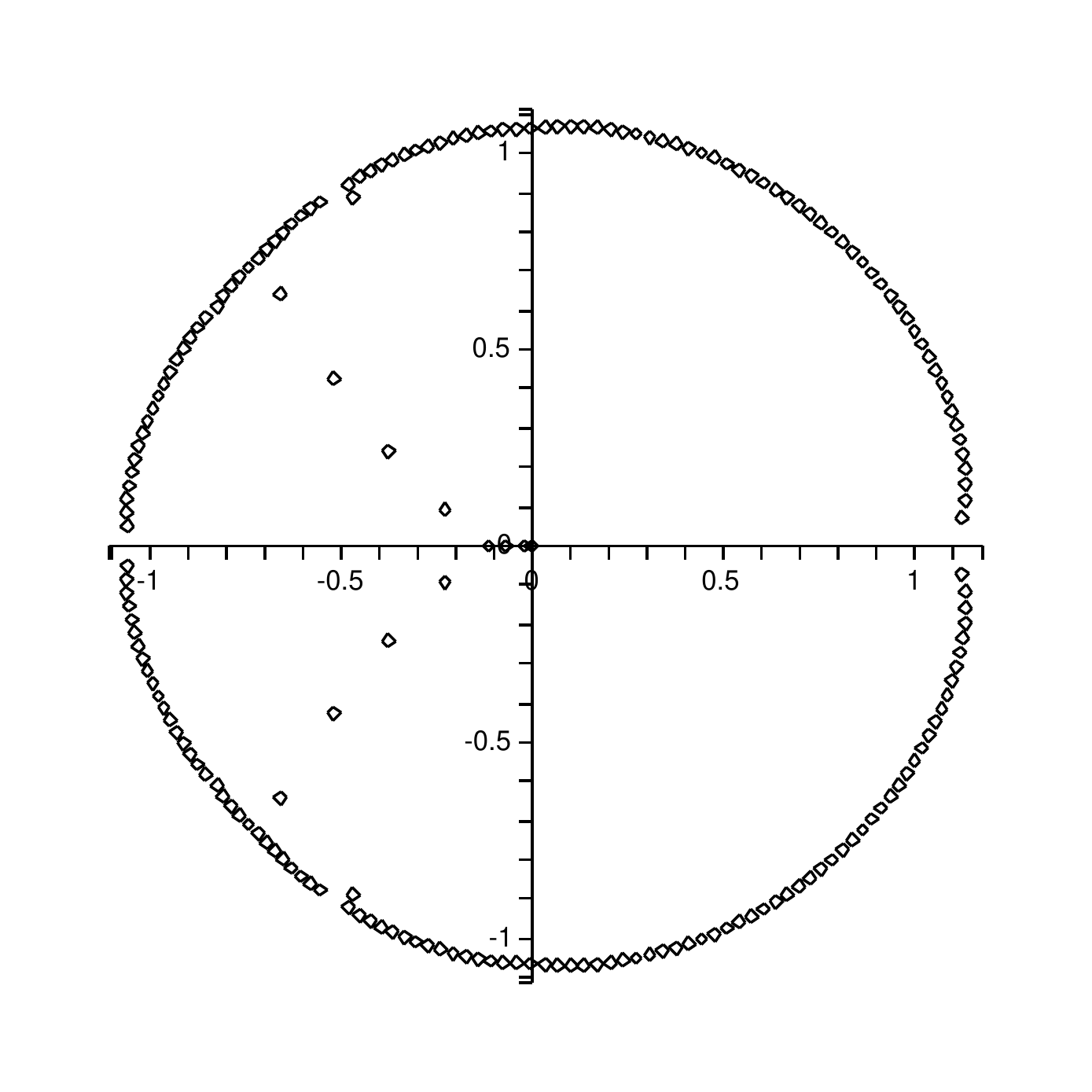} 
\qquad
\includegraphics[height=6.5cm,width=6.5cm]{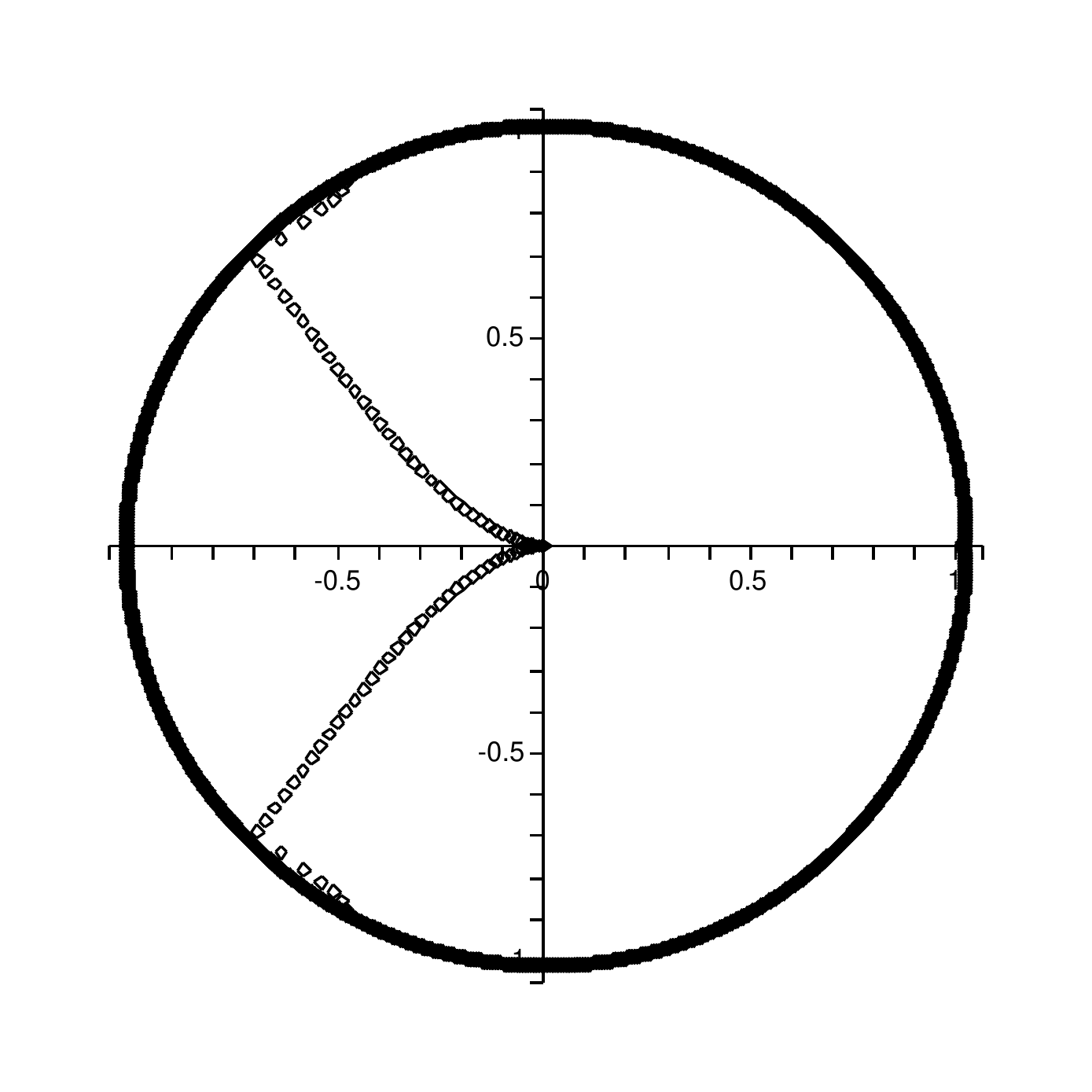}
 \caption{(a) Partition Polynomial Zeros for Degree 200; \quad (b) for Degree 10,000}
 \end{center}\end{figure}

Since the zeros of $F_n$ are symmetric about the real axis,
many times we restricted our attention to the upper half-plane.

We single out many intriguing features of this plot:
\begin{enumerate}
\item
The zeros are clustering about the unit circle;
\item
There is a sparse of zeros in the left half plane;
\item
There are gaps near the angle $2 \pi/3$ and at $x=1$ and $x=-1$.
\item
There may be two scales of zeros for those either outside or inside the unit circle;
\item
There are $O(n)$ zeros outside the unit circle and $O(\sqrt{n})$ inside.
We use $\sqrt{n}$ rather than other powers of $n$ by  inspiration from the Hardy-Ramanujan
asymptotics from $p(n)$.
\item
Empirically we expect the order of convergence of the zeros to
be $O(1/\sqrt{ n})$  \cite{varga}.
\end{enumerate}

We compare these proposed features with the computed zeros of
$F_{10000}$ in Figure 1(b).
At degree $10,000$  all three gaps now disappear while a second family
of zeros inside the unit disk appears. 
This is especially noteworthy at $x=1$ since $F_n(1)=p(n)$.
Furthermore, the order $O(\sqrt{n})$ of zeros is confirmed inside the disk.

These zeros were found using the 
MPSolve program described in \cite{bini} which
is our major tool  to obtain the zeros of high degree
polynomials. 
Its underlying algorithm is based on simultaneous approximation
of all the zeros and uses the Aberth iteration.
One  advantage of MPSolve is that
it handles integer coefficients with arbitrary large precision.
This ability is critical since the coefficients of the partition
polynomials  have hundreds of digits; for example, in Figure 2(b) there is a
plot of the number of digits of the 25,000 degree partition polynomial.
Note that the maximal coefficient 
 168 digits and that the polynomial is unimodal.
 For comparison in Figure 2(a), we included $F_{500}$ whose maximal number of
 digits is 19.

\begin{figure}[h!]
\begin{center}
\includegraphics[height=5cm,width=5cm]{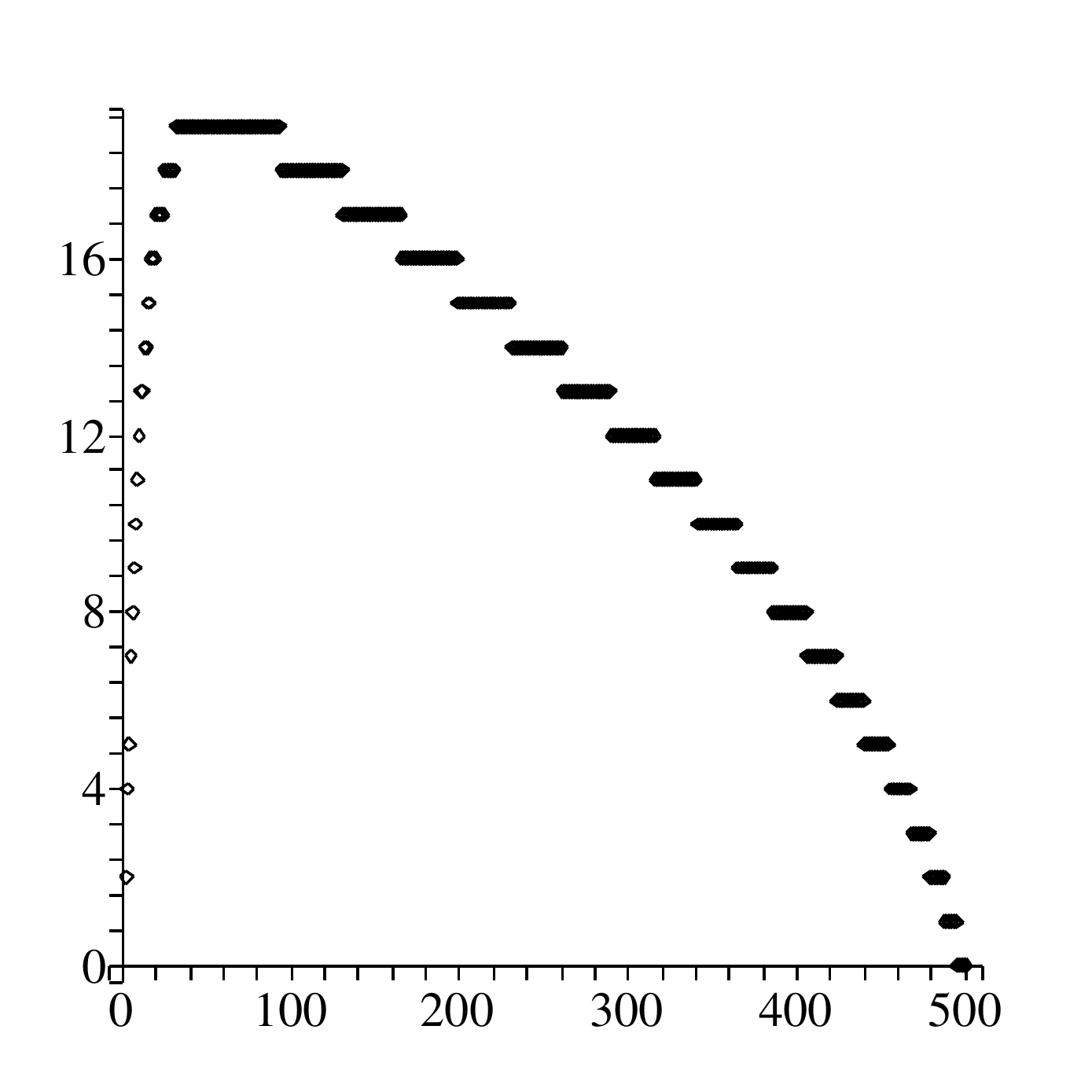}
\qquad
\includegraphics[height=5cm,width=5cm]{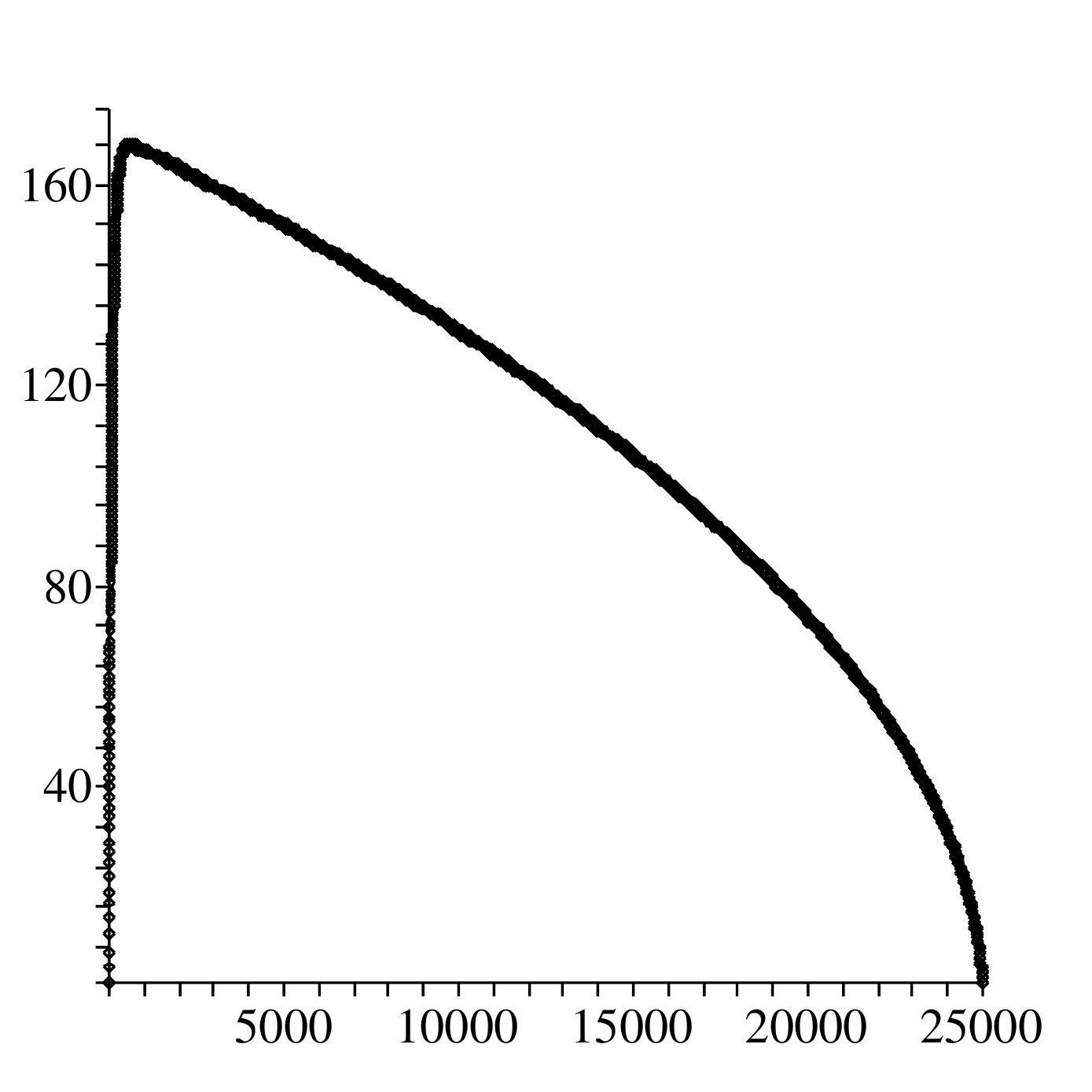}
\\
\caption{(a) Digits of the Partition Polynomials for Degree 500;\qquad
(b) for Degree 25,000 }
\end{center}\end{figure}

In working with such polynomials, a natural question is how
to test whether the zeros are relatively correct. The zeros obtained
have, say, 20 digits, yet the coefficients of the polynomial have
hundreds of integer coefficients. So there is no possibility of
ever checking directly that the zeros through evaluation. 
On the other hand, the zeros do pass  a simple numerical test.
We observe that the sum of the zeros equals the negative
of the coefficient of $x^{n-1}$ in $F_n(x)$ which is always 1.

The sum of the computed zeros are:
at degree 6000,  $-1.0000000000000000014085$;
at degree 10,000, $-0.9999999999999999860563$;
at degree 15,000, $- 1.00000000000000002067132+ 2 \times 10^{-19}i$;
and
at degree 20,000, $-1.00000000000000002323212$.
For $F_n$, the coefficient of $x^{n-2}$ is $2$ and equals
 the second elementary symmetric function of the zeros.
At degree 6000, for the computed zeros this evaluates to
 $1.99999999999999999377164 + 5.62 \times 10^{-26} i$.

Furthermore,
the zeros that MPSolve produces coincide extremely well with the
asymptotics we discovered so we have much confidence in the
computation.
In our work there was a complementary interplay between
developing the asymptotics and numerically determining their
limiting behavior and their densities.

The appealing question of determining the limiting behavior
of the zeros is intimately connected with the analytic problem of
finding asymptotic formulas for these polynomials as their
degrees go to infinity.

\section{Introduction to the Zero Attractor}\label{intro:zero}

We now formalize what we mean by the convergence of the
zeros.
Let ${\mathcal Z}(F_n)$ denote the finite set of zeros of
the polynomial $F_n$. Then the {\sl zero attractor} $\mathcal A$
of the polynomial sequence $\{F_n\}$ is the limit of
${\mathcal Z}(F_n)$ in the Hausdorff metric on the compact
subsets of ${\mathbb C}$.  
In Figure 3(a), there is the plot of the full zero attractor as well as a closeup of
the upper left-hand quarter of the attractor in Figure 3(b). We will describe the
curves in these plots in Section \ref{section:explicit_zero}.

\begin{figure}[h!]
\begin{center}
\includegraphics[height=6cm,width=6cm]{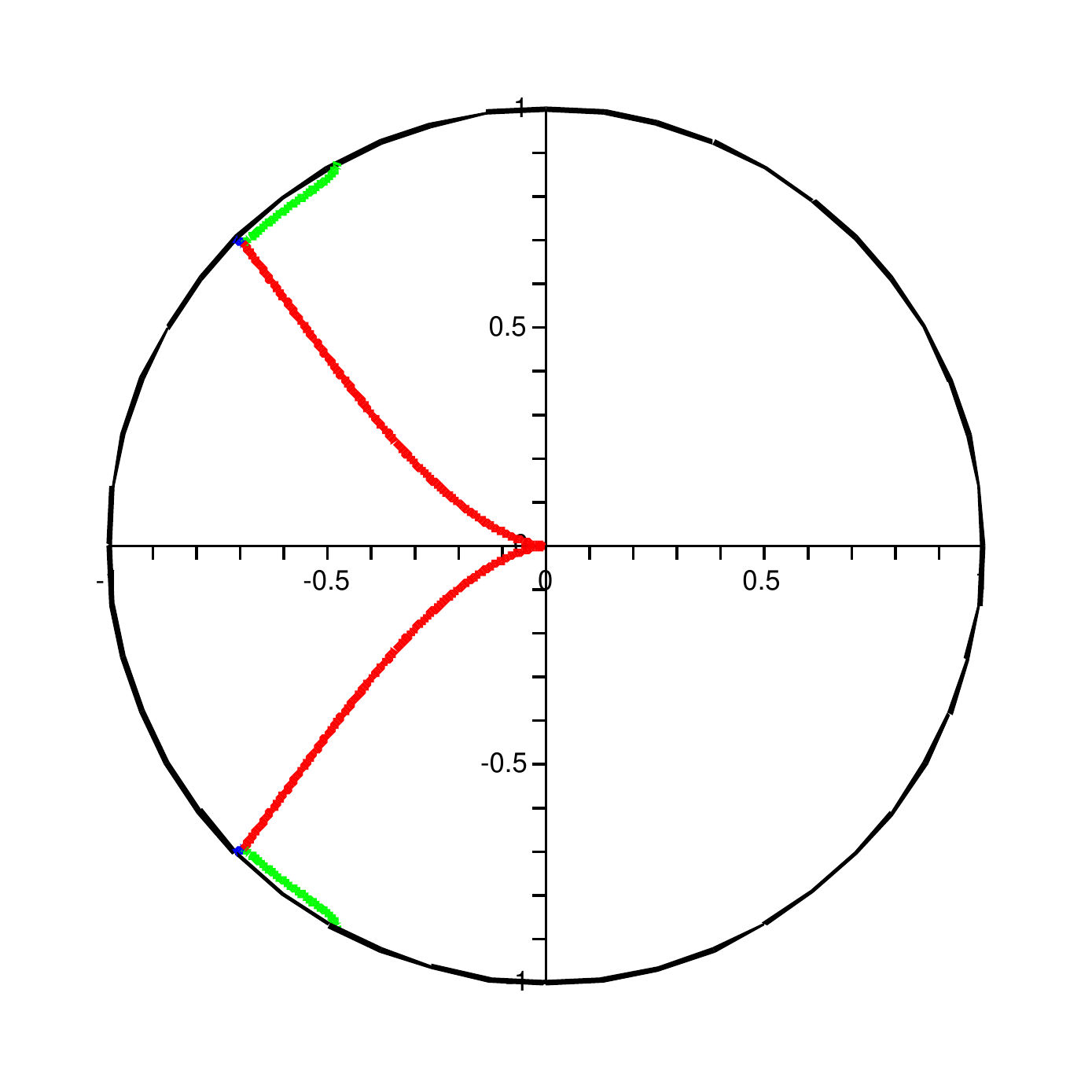}
\qquad
\includegraphics[height=6cm,width=6cm]{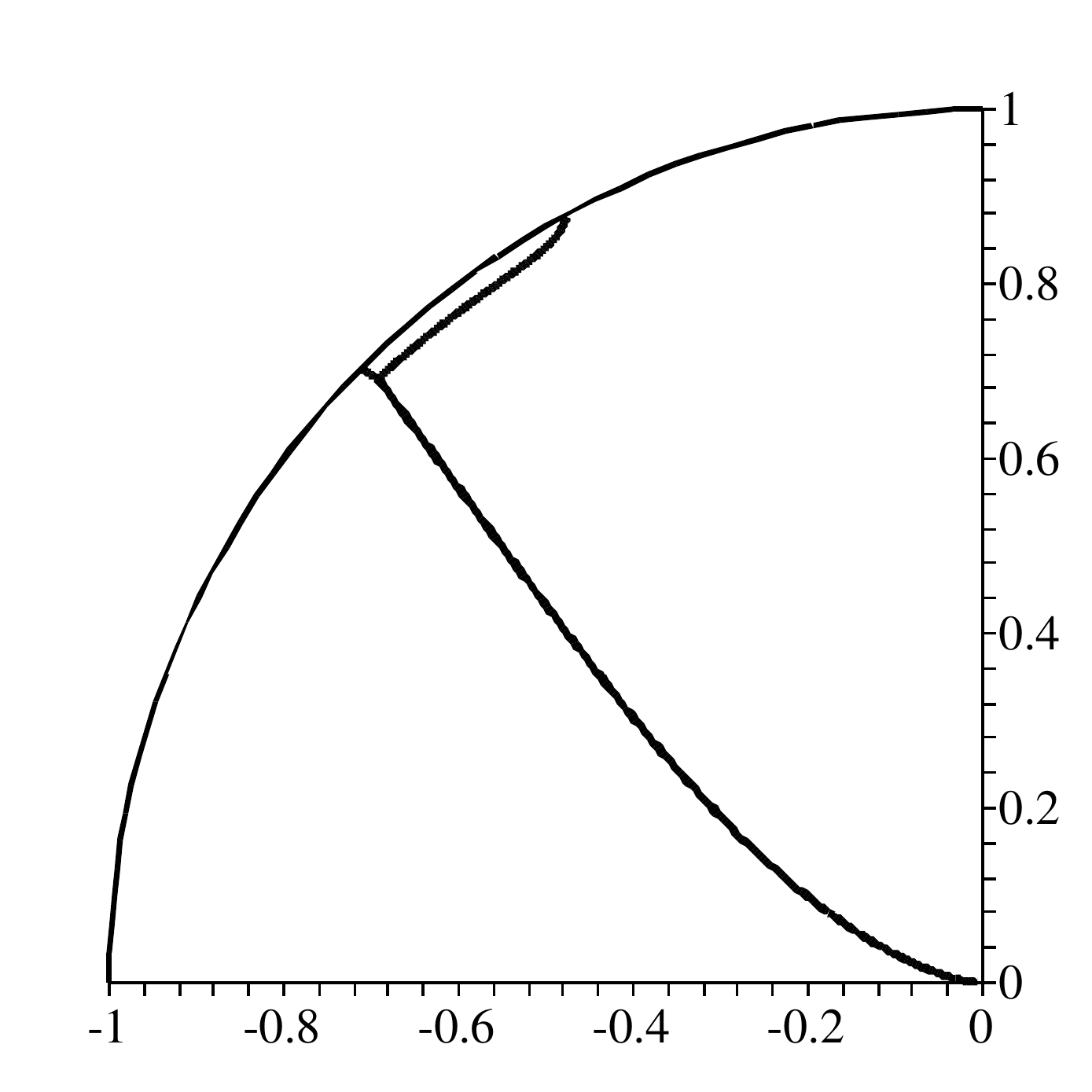}
\caption{
(a) Full Zero Attractor; \qquad (b) Zero Attractor in Left Half Plane }
    \end{center}\end{figure}

  It may be  hard to see but there are actually three curves that make up the zero attractor in the upper half plane.  Furthermore
  the common intersection point of these curves we call the {\sl triple point} $T$.
 Its  polar coordinates are $r=0.97833 70882,
\theta=2.356797156$
and its rectangular coordinates are $x=-0.6922055811$
and $y=0.6913717463$.

  The zeros converge very slowly to this special point.
 In Figure 4, there are plots of the zeros near the triple point $T$ for 
 degrees $400$, $5,000$, and $50,000$:

  \begin{figure}[h!]
   \begin{center}
\includegraphics[height=5.2cm,width=5.2cm]{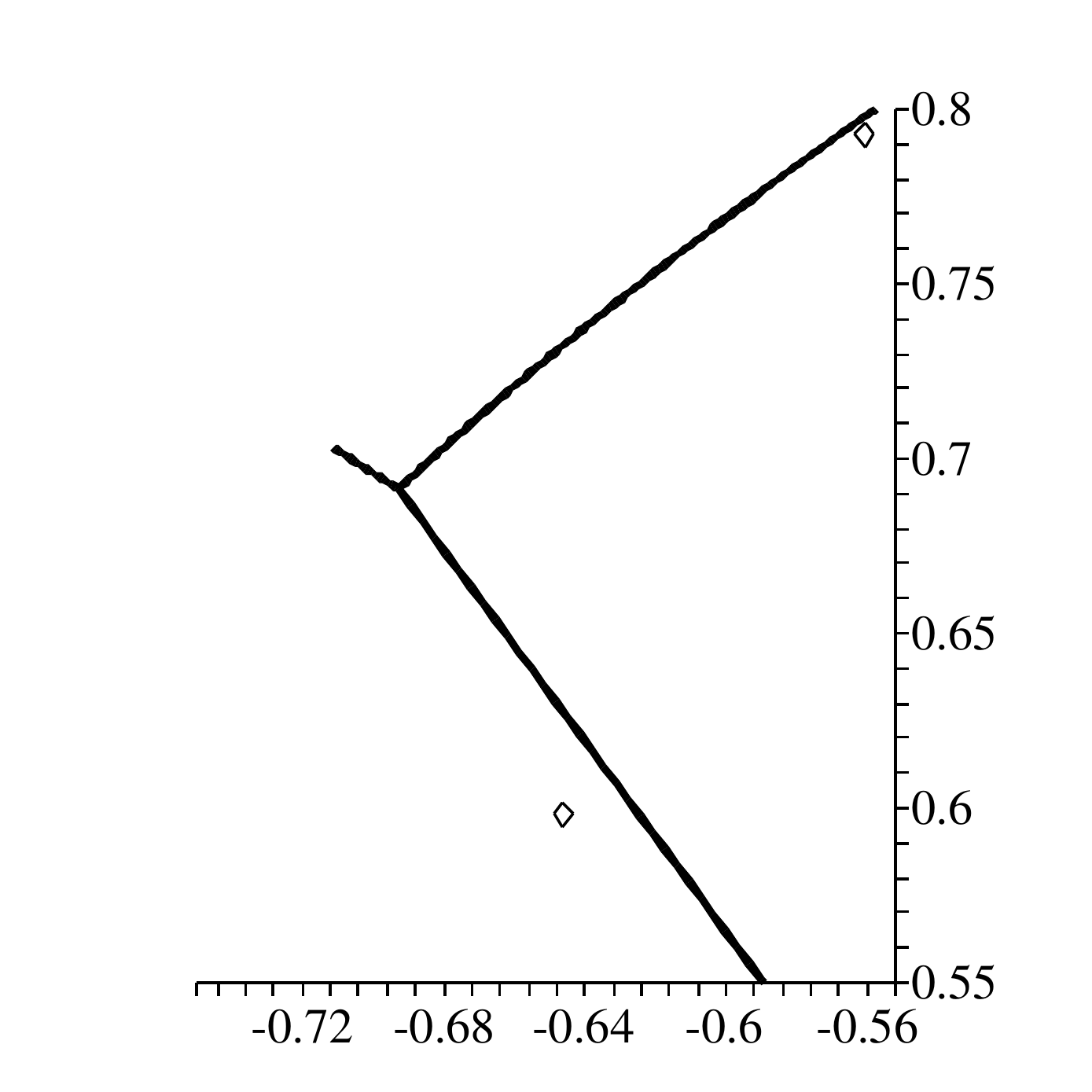}
\
\includegraphics[height=5.2cm,width=5.2cm]{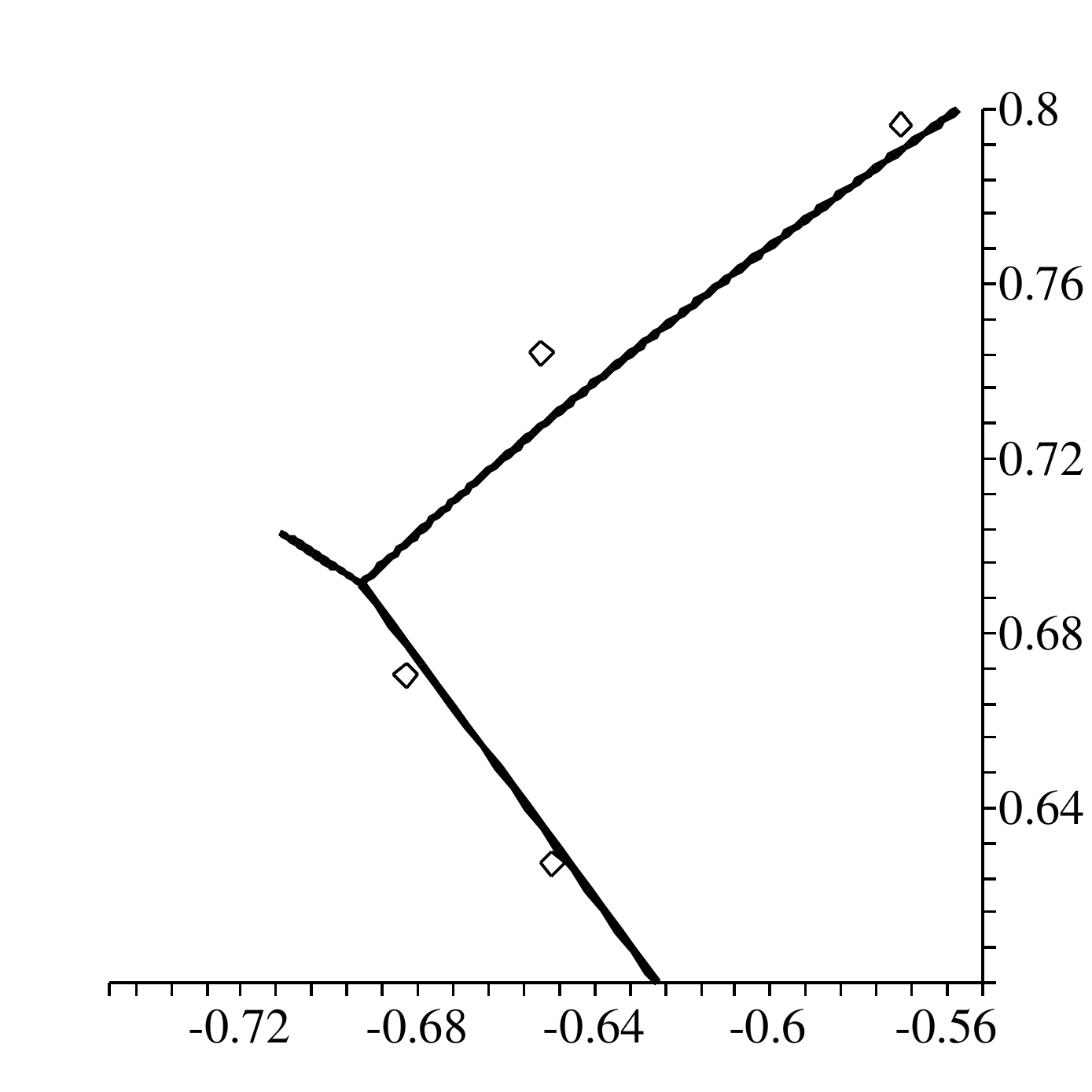}
\
\includegraphics[height=5.2cm,width=5.2cm]{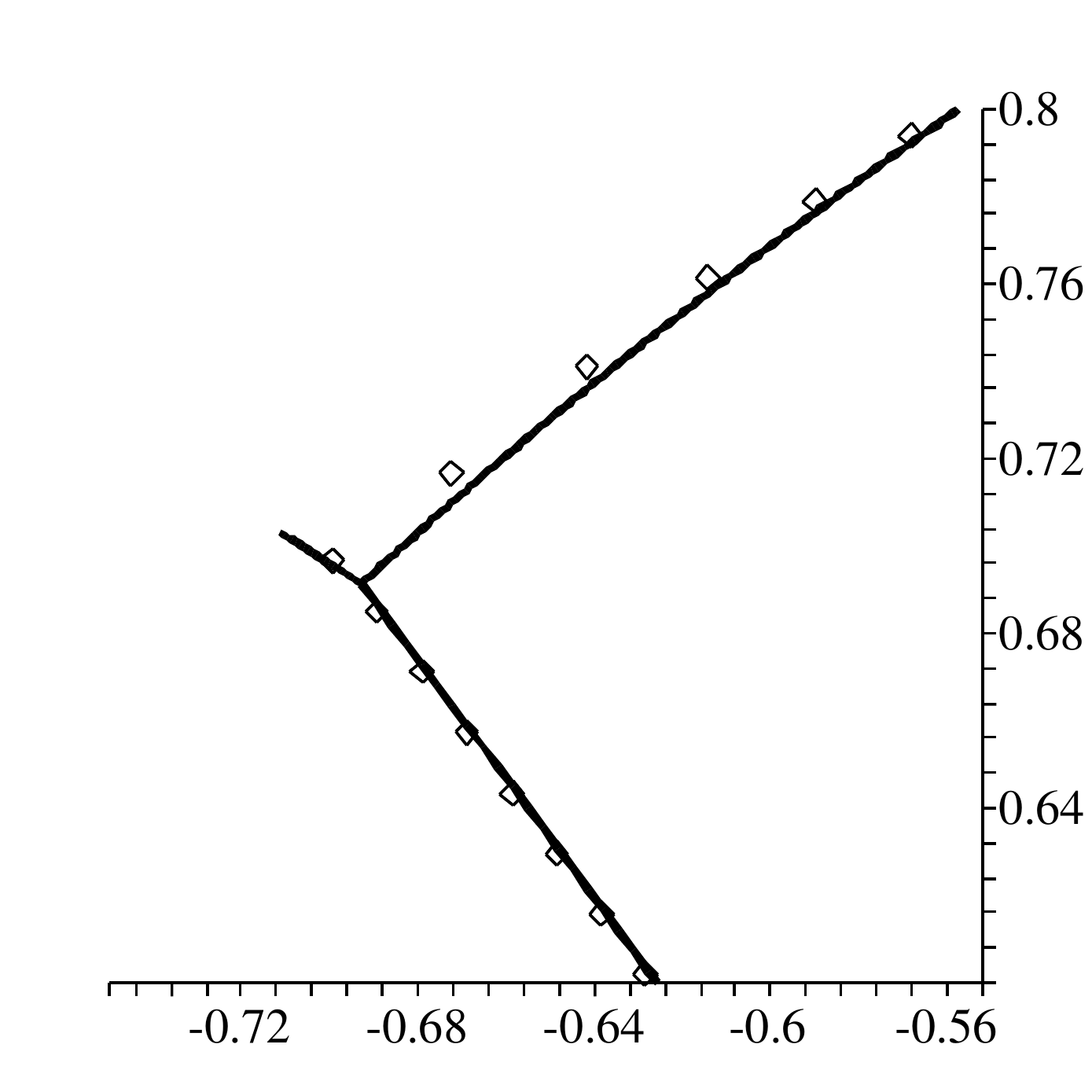}
\caption
{
Zeros for Degrees 400, 5,000, and 50,000 Near the Triple Point $T$
}
\end{center}\end{figure}

In Table 1, we give the total number of all zeros for $F_n$
inside the unit disk which empirically confirms that their order
inside the disk is $O(\sqrt{n})$. 
Its first entry is the count of all zeros inside the unit disk.
The second entry is the count of all zeros that lie in
 $Q_2$, the closed quarter
unit in the second quadrant excluding the unit circle.
The next entries are
 the counts  for the three families of zeros 
 near the three curves seen in Figure 3 that
 in the second quadrant.
These families will be identified in Section  \ref{section:explicit_zero}; 
 for now, we call
 them either Family 1, 2, or 3. Note that Family 1 includes zeros
 that lie along the real axis.

\begin{table}
\caption{Number of Zeros Inside the Unit Disk}
\label{page:table}
\begin{center}
\begin{tabular}{| r || c|c  |c |c| c|| r  |}
\hline
{\bf Degree} & $\#$ {\bf Zeros} &{\bf  All Zeros in} $Q_2$ & {\bf Family} 1 & {\bf Family} 2 & {\bf Family} 3 
& {\bf Prediction} \\
\hline
5000 & 64 & 36 & 32  & 4 & 0 & 64.8\\
10000& 92 & 51& 45  & 6 & 0 & 91.7 \\
15000& 112 & 61 &53  & 7 & 1 & 112.2 \\
20000&130 & 71 & 62 & 8 & 1 &129.6 \\
25000&146 & 79 & 69 & 9 & 1 & 144.9 \\
30000&160 & 87 & 76 & 10 & 1 & 158.7 \\
35000&172 & 93 & 81 & 11 & 1 &171.5\\
40000&184 & 99 & 86 & 12 & 1 & 183.3\\
50000&204 & 109 & 95 & 13 & 1 & 204.9 \\
60000&226 & 121 & 104 & 15 &2 & 224.5 \\
70000&242&  129 & 112&16&1 & 242.5\\
\hline
\end{tabular}
\end{center}
\end{table}

The predicted number of zeros in Table 1 is given by the empirical formula
\begin{equation}\label{eq:least_squares_law}
\# \textrm {  Zeros inside the unit disk} \sim 0.9154 \sqrt{n}
\end{equation}
obtained by fitting a multiple of the power law $\sqrt{n}$ to
the zeros to the first four polynomials at  degrees $5000, 10,000$, $15,000$, and $20,000$. Another estimate for these zeros is given in
equation (\ref{eq:square_root_constant}).

\section{Potential Theory, Statistical Mechanics,
and Limits of Zeros}\label{section:potential}
 
 If $\mu$ is a finite measure with compact support on
 ${\mathbb C}$, then its {\sl logarithmic potential function}
 $U^\mu(z)$ is given as
 \[
 U^\mu(z)= \int_{\mathbb C} \log \frac{ 1}{ |z-t|} \, d\mu(t).
 \]
 The connection with  polynomial zeros comes from the 
 fact that if $a_1, \cdots, a_n \in {\mathbb C}$ and
 $\mu = \frac{1}{n} \sum_{k=1}^n \delta_{a_k}$, then
 $U^\mu(z) = \frac{1}{n} \ln |P(z)|$, where $P(z)=
 (z-a_1) \cdots (z-a_n)$.
 
 The hope is 
 that the  logarithmic potentials
 for a sequence of polynomials 
 will converge to the density measure
 on the zero attractor.
 Classically, one studies  the first limit below while we need the following
 second limit below as well because the number of zeros inside the unit disk
 is $O(\sqrt{n})$:
 \begin{equation}\label{eq:ln_limit}
 \lim_{n \to \infty} \frac{\ln[ F_n(x)]}{n},
 \quad
  \lim_{n \to \infty} \frac{\ln[ F_n(x)]}{\sqrt{n}}.
 \end{equation}

 From the Yang-Lee theory of phase transitions in statistical
 mechanics, the limit of ${\ln[ F_n(x)]}/{n}$ gives the complex
 free energy and its zeros indicate the presence of phase transitions.
 
 We found Alan Sokal's formulation  \cite{sokal} to be very helpful
because it  allows more general normalizations than $n$.
 Here is his result:
 
 \begin{theorem}$\,$
 [Sokal]
 Let $D$ be a domain in $\mathbb C$, and let $z_0\in D$.
 Let $\{g_n\}$ be a sequence of analytic functions in $D$,
 and let $\{a_n\}$ be a sequence of positive reals such that
 $\{ | g_n |^{a_n}\}$ are uniformly bounded on the compact 
 subsets of $D$.
 Suppose there does not exist a neighborhood $V$ of $z_0$
 and a function $v$ on $V$ that is either harmonic or else
 identically $- \infty$ such that
 \[
 \liminf_{n\to \infty} a_n \ln | g_n(x)|
 \leq v(x)
 \limsup_{n\to \infty} a_n \ln | g_n(x)|
 \]
 for all $z\in V$. Then $z_0$ lies in $\liminf {\mathcal Z}(g_n)$; 
 that is, for all $n$ sufficiently large, there exists zeros $z_n^*$ of
 $g_n$ with $\lim_{n\to \infty} z_n^* = z_0$.
 \end{theorem}
 
 Below we will state what  the limits in equation (\ref{eq:ln_limit})
 are and see that they determine on three regions inside the disk.
This theorem tells us that the zeros accumulate on the boundaries
of these regions. Note: these limits are very difficult to compute.

 \section{Explicit Description of Zero Attractor}\label{section:explicit_zero}
 
 Let $\Li_2(x)$ denote the dilogarithm of $x$; that is,
 $\Li_2(x) = \sum_{n=1}^\infty \frac{x^n}{n^2}
 =
 - \int_0^x \frac{ \ln(1-t)}{t} \, dt$.
 We introduce a related family of (sub)harmonic functions
 \begin{eqnarray*}
 f_1(x) = \Re \left[ \sqrt{ \Li_2(x)} \right],
 \quad
 f_k(x) = \frac{1}{k} \Re \left[ \sqrt{ \Li_2(x^k)} \right]
 \end{eqnarray*}
 which are harmonic in different sectors of the unit disk.
These functions play a crucial role in the asymptotics
of the partition polynomials.
Next we introduce the curves $C_{k, \ell}$ inside the unit disk
where $f_k(x) = f_\ell(x)$
which we will informally call $\sqrt{\Li_2}$-curves.
The zero attractor consists of portions
of $C_{12}$, $C_{13}$, and $C_{23}$.

We can now describe how the plots of the zero attractor are done.
Since the $\sqrt{\Li_2}$-curves $C_{k,\ell}$ are level curves of the real part
of a harmonic function,  these curves are naturally solutions to
an initial value problem. For convenience, write $L_k(z) =
\sqrt{\Li_2(z^k)}/k$ and  $L_k(z)'$ is just its usual derivative. Then
\[
\frac{dy}{dx} = \frac{ \Re [  L_k(x+iy)' -   L_\ell(x+iy)' ]  }
{\Im[   L_k(x+iy)' -   L_\ell(x+iy)'   ]  }, \quad y(x_0)=y_0.
\]
The initial condition needs to be found numerically; typically, by a value
on the unit circle where $\Li_2(e^{it})=u(t)+iv(t)$ is known in closed form:
\begin{equation*}
u(t)=\sum_{n=1}^{\infty }
\frac{\cos n t }{n^{2}}=\frac{3 t^{2}-6 t \pi +2\pi ^{2}}{12},
\quad
v(\theta )=\sum_{n=1}^{\infty }\frac{\sin n t }{n^{2}}
=-\int_{0}^{ t }\ln (2\sin \frac{\xi }{2})\, d \xi , \quad 0 \leq t \leq 2\pi .
\end{equation*}
Since the $\sqrt{\Li_2}$-curves are very smooth and even Euler's method gives an
useful plot for them.

In establishing the role of
the $\sqrt{\Li_2}$-curves 
for  the zero attractor, though, it is
 better to treat them as boundaries for regions
of majorization among the functions $f_k$
inside the upper unit disk. We define
\begin{equation}\label{eq:regions}
{\mathcal R}(1) = \{ x  : f_1(x) \geq f_2(x), f_3(x)\},
{\mathcal R}(2) = \{ x  : f_2(x) \geq f_1(x), f_3(x)\},
{\mathcal R}(3) = \{ x  : f_3(x) \geq f_1(x), f_2(x)\}.
\end{equation}
Then $\partial {\mathcal R}(1) $ consists of portions of the
unit circle, $C_{12}$, $C_{13}$, and $[0,1]$; 
$\partial {\mathcal R}(2)$ consists of portions of the unit
circle, $C_{12}$, $C_{23}$, and $[-1,0]$; 
$\partial {\mathcal R}(3)$ consists of portions of the unit
circle, $C_{13}$, and $C_{23}$.

\begin{figure}[h!]
\begin{center}
 \includegraphics[height=7cm,width=7cm]{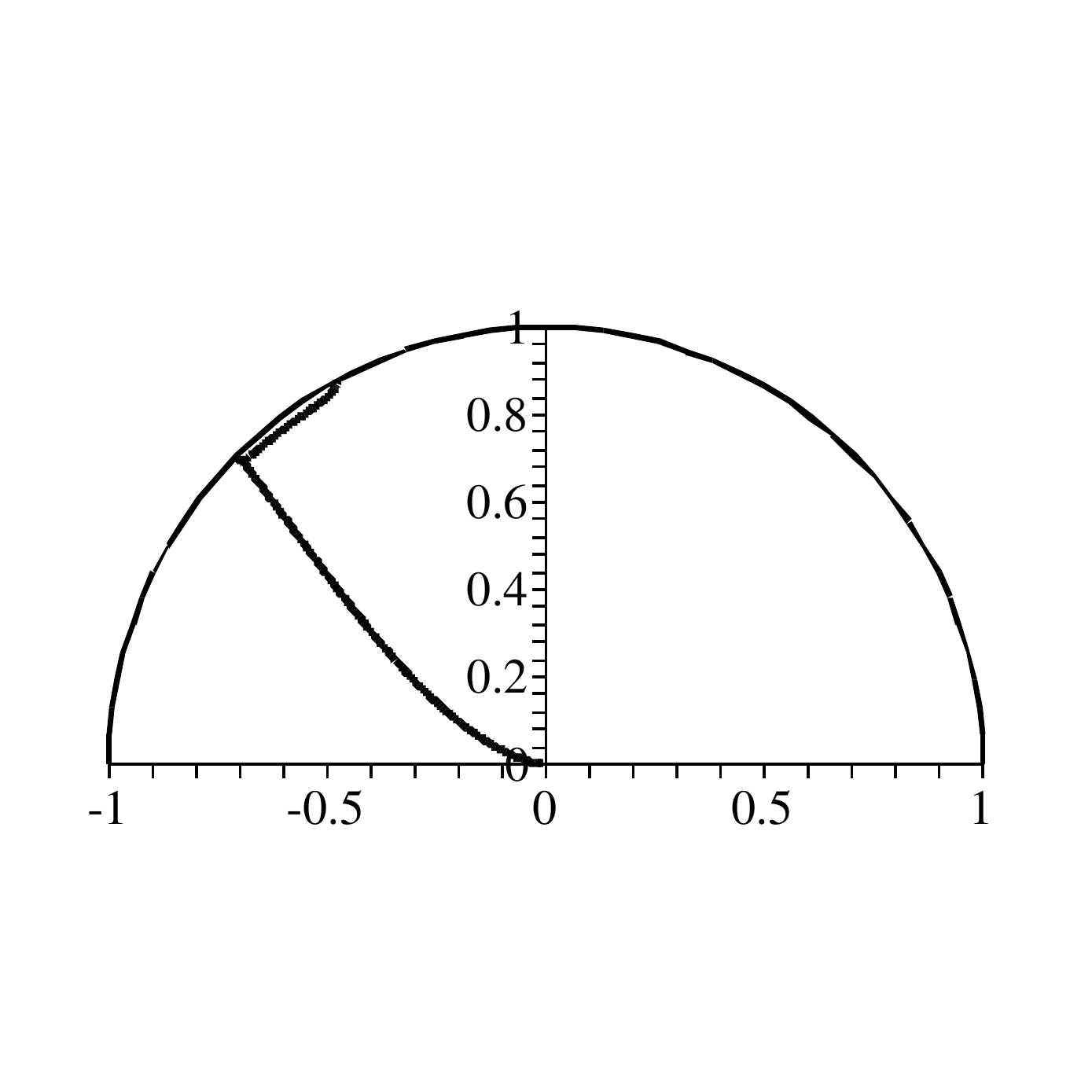} \qquad\qquad
  \includegraphics[scale=0.35]{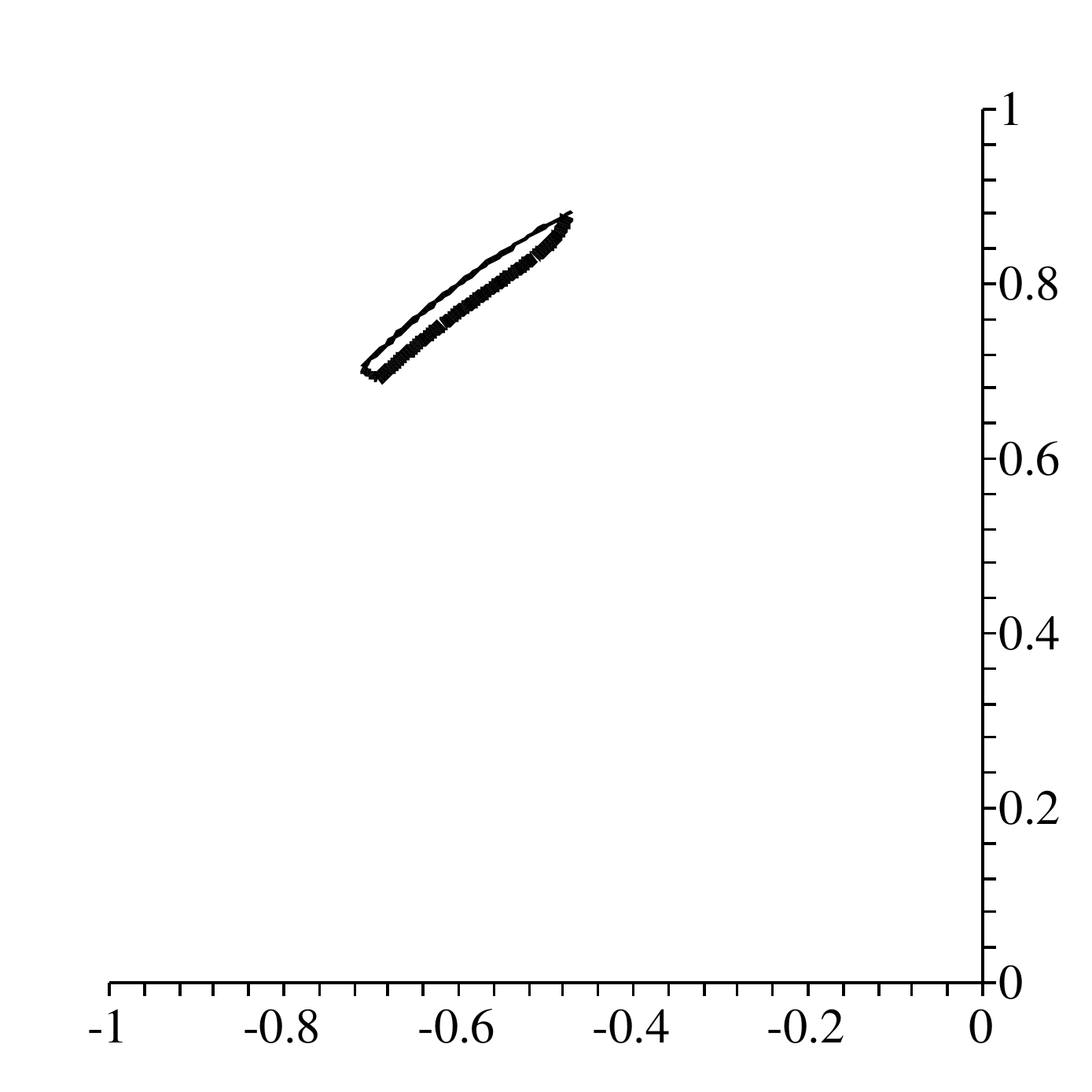} 
   \caption{(a) The Three Regions ${\mathcal R}(1)$,
  ${\mathcal R}(2)$, ${\mathcal R}(3)$; \hfill
  (b) Region ${\mathcal R}(3)$ Separately}
   \end{center}\end{figure}

  The realization that the zeros accumulation on the boundary curves
  made us examine the majorization among $f_1$, $f_2$, and $f_3$ on
  the unit circle. Since these are known in closed form, we found 
  that $f_1(e^{it})$ dominates on the arc $t\in [0,\theta_{13}]$, $f_3(e^{it})$ dominates
  on $t\in [\theta_{13}, \theta_{12}]$, and $f_2(e^{it})$ dominates on $t\in [\theta_{12}, \pi]$
  where
  $\theta_{13}= 2.066729664 < 2 \pi/3=2.094395103$,
\
$\theta_{12}=2.2536266 < 3 \pi/4=2.356194490$,
and
$\theta_{23}=2.361704176 > 3\pi/4$.
The following Theorem together with Sokal's result guarantees 
that there are no further zero families.

  \begin{theorem}\label{theorem:domination}
  (a) On the unit disk,
  $f_k(x) \leq \max[ f_1(x),f_2(x),f_3(x)]$, for all $k \geq 4$.
  \\
  (b) For $x\in {\mathcal R}(1)$, $f_1(x) \geq f_k(x)$, for all $k \geq 2$.
  \\
  (c) For $x \in {\mathcal R}(2)$, $f_2(x) \geq f_k(x)$, for all $k\geq 3$
  and $k=1$.
  \\
  (d) For $x \in {\mathcal R}(3)$, $f_3(x) \geq f_k(x)$, for all $k \geq 4$
  and for $k=1,2$.
  \end{theorem}
  
  The proof is intricate and lengthy and requires subharmonic function theory and
  special properties of conformal maps.  The theorem, though,  is easy to verify numerically
 in special cases.

  We  next record the values of the normalized limits
  of $\ln |F_n(z)|$. 
  
  \begin{theorem}
  Outside  the unit disk,
  $\displaystyle
  \lim_{n \to \infty}
   \frac{ \ln [ F_n(z)]}{n} = \ln z,\,
\lim_{n \to \infty}
 \frac{ \ln | F_n(z) | }{n} = \ln |z|.
  $
  \\
  Inside the unit disk,
  $\displaystyle
  \lim_{n \to \infty} 
  \frac{ \ln [ F_n(z)]}{n} = 0
  $
  and
  $\displaystyle
\lim_{n \to \infty}
 \frac{ \ln | F_n(z) | }{n} =0.
  $
  \end{theorem}
  
  By Sokal's theorem, we now know that the $O(n)$ contribution
  to the zero attractor is the unit circle.  The delicate asymptotic
  expansions needed to establish the following theorem are outlined
  in Section \ref{section:asymptotics_inside}.
  
  \begin{theorem}\label{theorem:density}
   
  \begin{enumerate}
  \item
  On region ${\mathcal R}(1)$,
$\displaystyle
\lim_{n \to \infty}
 \frac{ \ln [F_n(z)]}{ 2\sqrt{n}}=  \sqrt{\Li_2(z)}$ and
$\displaystyle
\lim_{n \to \infty}
 \frac{ \ln |F_n(z)|}{ 2\sqrt{n}}=  \Re  \sqrt{\Li_2(z)}.
$
\item
On region ${\mathcal R}(2)$,
$\displaystyle
\lim_{n \to \infty}
 \frac{ \ln [F_n(z)]}{ 2\sqrt{n}}= \tfrac{1}{2}  \sqrt{\Li_2(z^2)}$ and
$\displaystyle
\lim_{n \to \infty}
 \frac{ \ln |F_n(z)|}{ 2\sqrt{n}}= \tfrac{1}{2} \Re  \sqrt{\Li_2(z^2)}$.
\item
On region ${\mathcal R}(3)$,
$\displaystyle
\lim_{n \to \infty}
 \frac{ \ln [F_n(z)]}{ 2\sqrt{n}}= \tfrac{1}{3}  \sqrt{\Li_2(z^3)}$ and
$\displaystyle
\lim_{n \to \infty}
 \frac{ \ln |F_n(z)|}{ 2\sqrt{n}}= \tfrac{1}{3} \Re  \sqrt{\Li_2(z^3)}.
$
\end{enumerate}
  \end{theorem}

The above theorem shows that while the limit 
$ \frac{ \ln |F_n(z)|}{ 2\sqrt{n}}$ is continuous on the unit disk
the derivative of the limit fails to be continuous exactly on the
boundaries of the three regions. Hence we have determined
the zero attractor inside the unit disk:

\begin{theorem}
The zero attractor for the partition polynomials consists exactly
of the unit circle together with the boundaries of the three regions
${\mathcal R}(1)$, ${\mathcal R}(2)$, and ${\mathcal R}(3)$
inside the unit disk.
\end{theorem}

\section{Zero Density}\label{section:zero_density}

The zero densities are determined by our following theorem:

\begin{theorem}
{\bf Density Theorem}
Let $G$ be a conformal map from a neighborhood of an analytic arc
$C$ to a neighborhood sector $S$ of an arc of the unit circle 
such that $G(C)=S$.
Let $\{T_n(x)\}$ be a sequence of analytic functions on $G^{-1}(S)$
such that
\[
T_n(x) = 1+ a_n(x) G^{c_n}(x)+e_n(x),
\]
where $\{c_n\}$ is an unbounded sequence of increasing positive
numbers and $\{a_n(x)\}$ is a sequence of analytic functions that
satisfy uniformly on $G^{-1}(S)$ as $n \to \infty$:
$|a_n(x)| \geq \delta > 0$, $\frac{a_n'(x)}{n a_n(x)} = o(1)$,
and $e_n(x) = o( a_n(x) G^{c_n}(x))$. Then for any $\epsilon>0$
all the zeros of $T_n(x)$ lie in $C_\epsilon$ for $n$ sufficiently large,
where $C_\epsilon$ is the $\epsilon$-neighborhood of the analytic arc $C$
and 
\[
\lim_{n \to \infty} \frac{1}{c_n}
\sum \{ \delta_{G(z)} : z \in {\mathcal Z}{(T_n)} \cap C_\epsilon \}
\to \mu
\]
where  $\mu$ is normalized Lebesgue measure on the unit circle restricted
to the circular arc $G(C)$. In particular, the zero
density measure $\nu_Z$ for the family $\{T_n\}$ is    $\nu_Z = \mu \circ G$
on the arc $C$.
\end{theorem}

To apply the Theorem to $\{F_n(x)\}$ requires special normalizations that
can be found in our paper \cite{boyer_goh_partition}.

\begin{table}
\caption{ Properties of the Boundary Curves $C_{12}\cap {\mathcal R}(1)$, 
$C_{13}\cap {\mathcal R}(1)$, $C_{23}\cap {\mathcal R}(2)$ }
\begin{center}
\begin{tabular}{| r || c |c |c|}
\hline 
{\bf Curve} &  $C_{12}\cap {\mathcal R}(1)$ & $C_{13}\cap {\mathcal R}(1)$ 
& $C_{23}\cap {\mathcal R}(2)$  \\
\hline\hline
{\bf Length} &  $0.9983742022$ & $0.2884481319$ & $0.02220012557$ 
\\
{\bf Density Mass} & $2.464527879$ & $0.367464849$ & $0.036529069$ 
\\
{\bf Map} $G_{k \ell}$ & $e^{ \sqrt{ \Li_2(x)} - \sqrt{\Li_2(x^2)}/2}$ & 
$e^{ \sqrt{ \Li_2(x)} - \sqrt{\Li_2(x^3)}/3}$ &
$e^{ \sqrt{ \Li_2(x^3)}/3 - \sqrt{\Li_2(x^2)}/2}$ 
\\
{\bf Circular Arc} & $[1, G_{12}(T)]$ & $[G_{13}(T), G_{13}(e^{i\theta_{13}})]$  
& $[G_{23}(T), G_{23}(e^{i\theta_{23}}   )]$
\\
{\bf Circular Arc} &  $[0,1.996701527]$ & $[1.388229082,1.755693930]$
& $ [ 1.077010447,1.113539516]$
\\
{\bf Relative Weight} &   $0.8591630301$  &   $0.1281025124$   &  $0.01273445753$
\\
\hline
\end{tabular}
\end{center}
\end{table}

Table 2 gives the essential information about the non-circular
boundary curves of the domination regions ${\mathcal R}(1)$,
${\mathcal R}(2)$, and ${\mathcal R}(3)$.  The entry ``length"
means the arc length of the curve while ``density mass" means
the total measure of the curve with respect to the zero density
measure $\nu_Z$.  By the Density Theorem, $\nu_Z$ is
the pull-back of Lebesgue measure on the unit circle under
the conformal mapping $G_{k \ell}$ where ${k \ell}= {12}$, ${13}$, or
${23}$.  The entry ``circular arc" means the image of $C_{k \ell}$
under $G_{k \ell}$.
Finally by ``relative weight", we mean
\[
w_{k, \ell} = 
\nu_Z(C_{k \ell})/ [ \nu_Z(C_{12}) + \nu_Z(C_{13}) + \nu_Z(C_{23})].
\]

We emphasize that the number of zeros inside the unit disk
is determined, of course, by the density measure and not by
the length  of the curves of the zero attractor. Here is
some simple evidence for this. For a particular polynomial $F_n$,
the number of its zeros near $C_{k \ell}$ is about $w_{k, \ell}$ 
out of the total number of zeros inside the unit disk in the
second quadrant.
We confirm this in Table 3. 
There should
be about $w_{13}+w_{23}$ or $14\%$ of zeros near curves $C_{23}$ and $C_{13}$.
At degree 70,000, there are 242  zeros inside the unit disk. About 14\% of them
should lie near these two curves or their reflections in the lower half plane.
This count yields 17 expected zeros which coincides with the computed number.

In Figure 6, there are plots  of  zeros along $C_{12}$ and $C_{13}$ of degree 50,000. 
On each plot the portion of the curve between adjacent boxes have equal
zero density measure. This can be understood graphically from counting
the number of zeros which are marked by a cross between boxes.
The horizontal axes in these plots are arclength. For $C_{12}$,
$s=0$ is the origin; for the other two, $s=0$ is the triple point $T$.

Figure 7 contains the density function for these curves.
Since the density function for $C_{12}$ becomes infinite at the origin,
we truncated its graph at $s=0.1$.
Furthermore, we see that the density functions for $C_{12}$ and $C_{13}$
are lightest  in a neighborhood of the triple point $P$ while heaviest near the origin
for $C_{12}$ and the unit circle at $e^{i \theta_{13}}$ for $C_{13}$. 
This explains the absence of
zeros near the triple point $T$ in Figure 5, even at degree 50,000.

\begin{figure}[h!]
\begin{center}
\includegraphics[height=5cm,width=5cm]{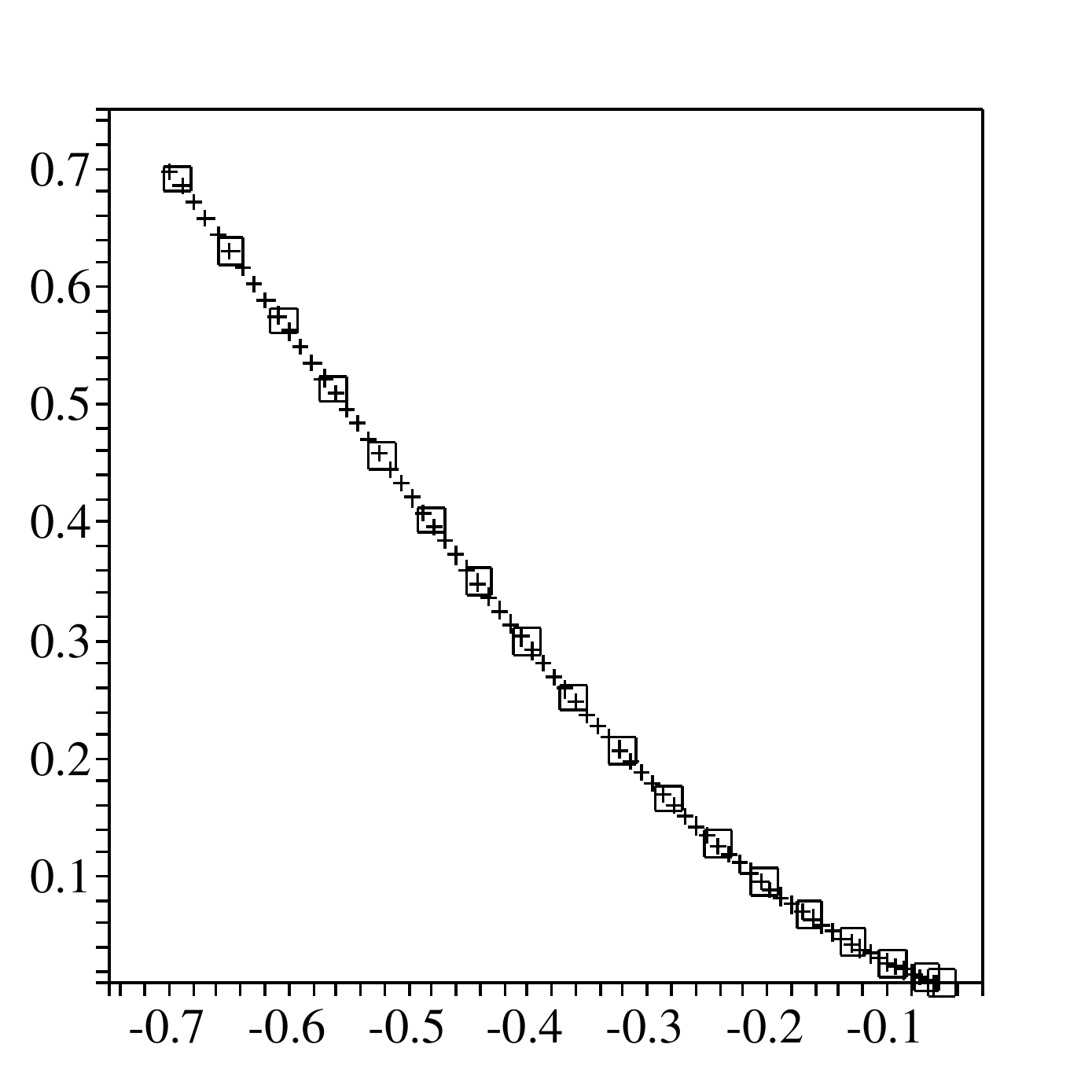} 
 \qquad
\includegraphics[height=5cm,width=5cm]{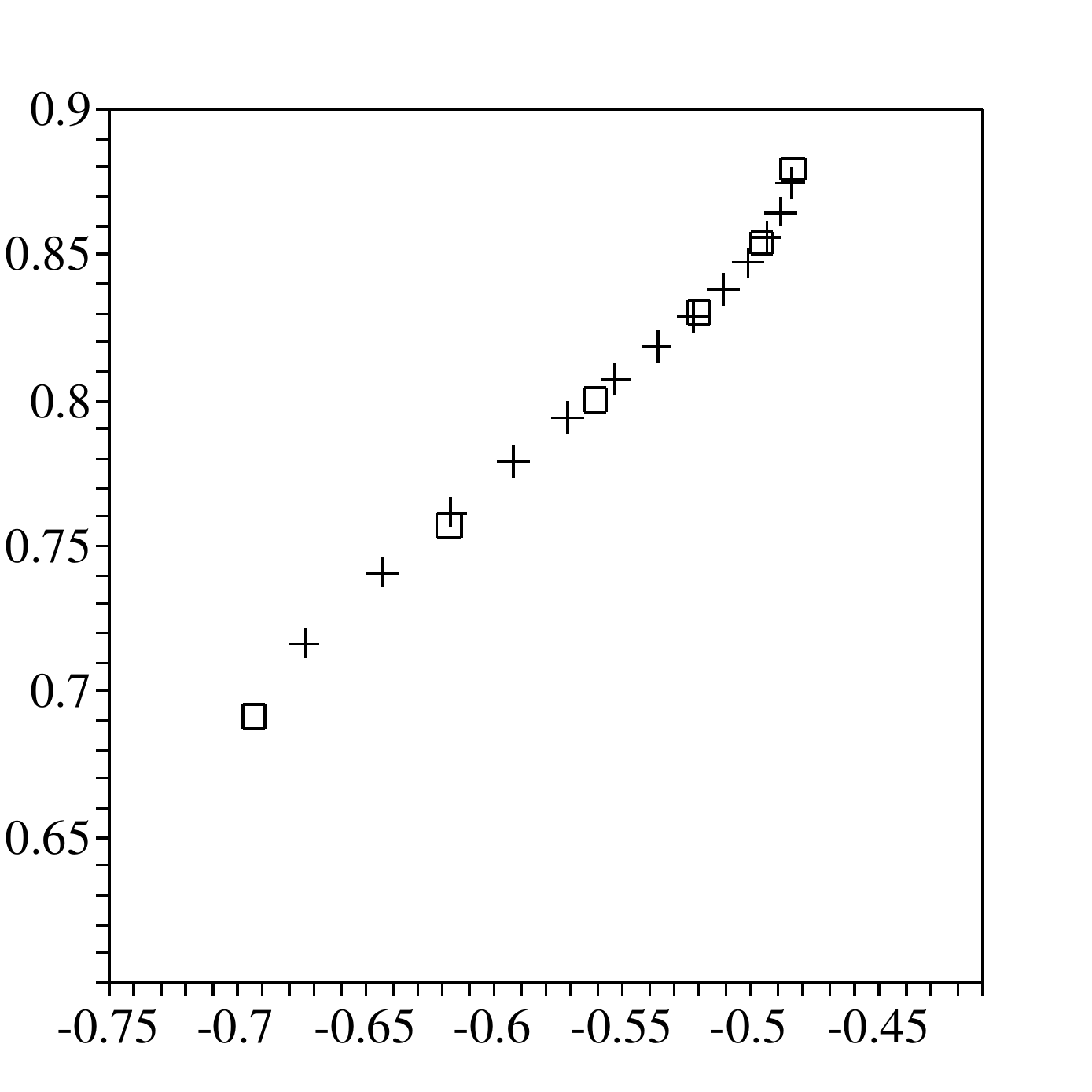}
 \caption
 {
 (a)  Zeros of Degree 50,000 Along Curve $C_{12}$ and (b) Along Curve $C_{13}$
  }
   \end{center}\end{figure}

  \begin{figure}[h!]
    \begin{center}
\includegraphics[height=4cm,width=4cm]{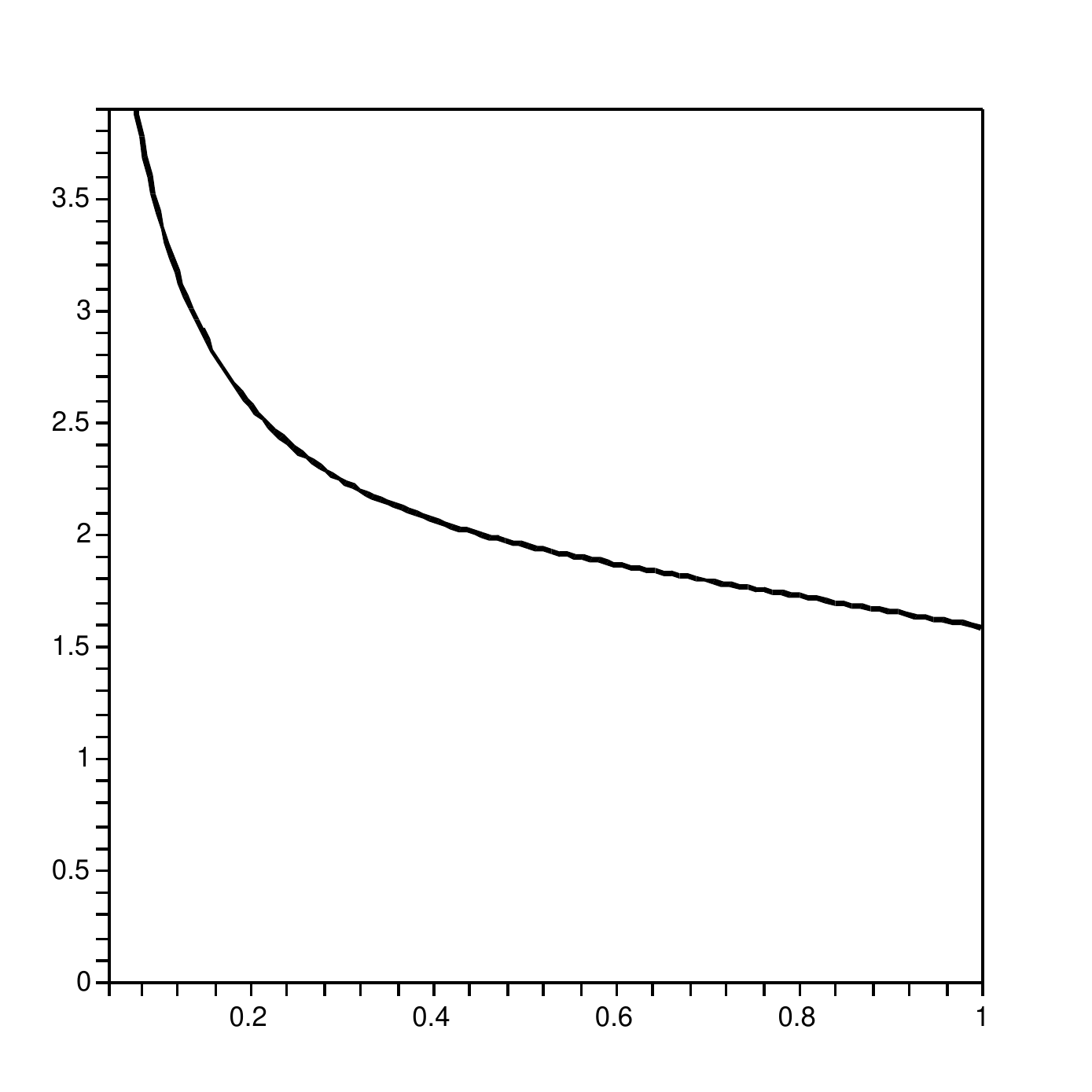} 
\
\includegraphics[height=4cm,width=4cm]{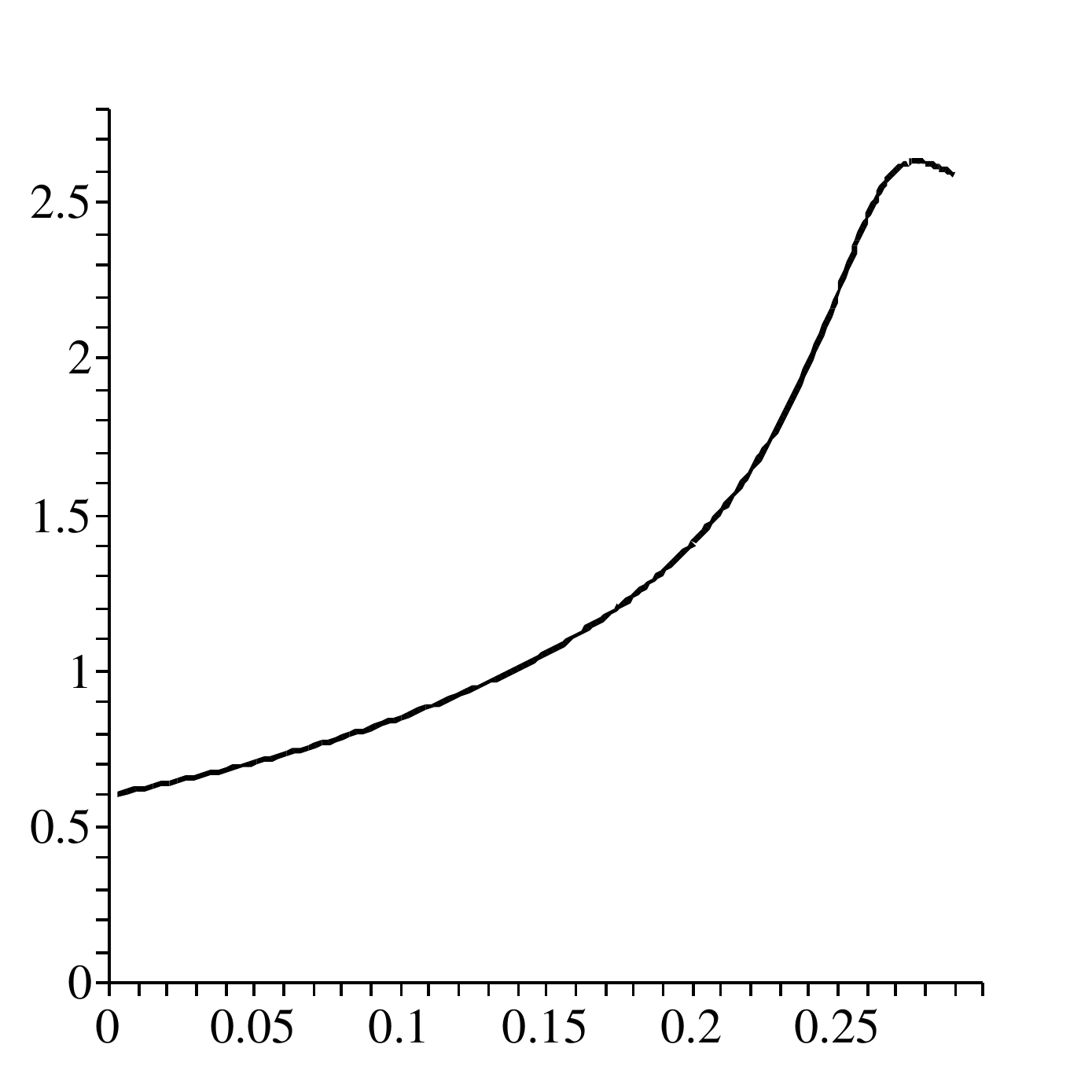} 
 \
\includegraphics[height=4cm,width=4cm]{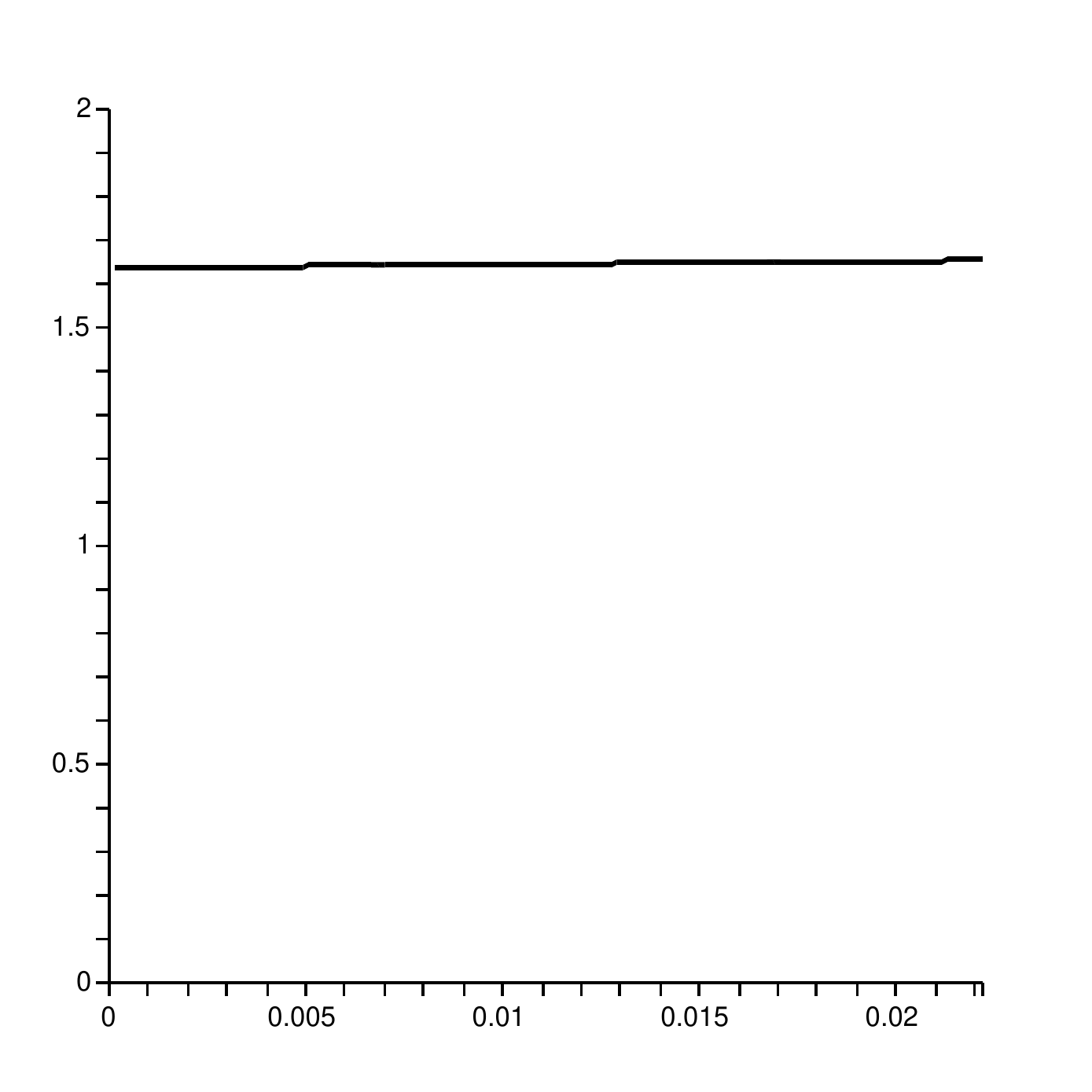} 
\caption
{
Density Functions as a Function of Arc Length for Curves
$C_{12}$, $C_{13}$, $C_{23}$
}
\end{center}\end{figure}

We can also explain why there are so few zeros for the $C_{23}$-family.
Its relative weight is $0.0127$. 
Now this makes two zeros near $C_{23}$ more likely for degree 70,000
over 60,000 yet the opposite occurs. At 70,000, the zero that does occur
is nearly at the center of the curve. A second zero would violate the near
uniform density. At 60,000, the two zeros that do occur are near the
two endpoints which is consistent with the zero density.
By equation (\ref{eq:least_squares_law}) and the relative weight of
$C_{23}$, we do not expect a third zero near this curve until
degree 190,000.

In Table 3,   the  predicted number of zeros near either of
the two $\sqrt{\Li_2}$-curves $C_{13}$ and $C_{23}$ is found  using 
the $\sqrt{n}$-approximation to the total number of zeros inside
the unit disk (see equation (\ref{eq:least_squares_law}))
 and the relative weights  of these two curves.
 
 Here is another way to get a nearly equivalent prediction of the
 number of zeros inside the unit disk:  $C \cdot \sqrt{n}$ where
 \begin{equation}\label{eq:square_root_constant}
 C = \frac{ \nu_Z(C_{12}) + \nu_Z(C_{13}) + \nu_Z(C_{23}) } { \pi}
 \simeq 0.9130788466.
 \end{equation}
 Note that $C$ is remarkably close to the value given by
 the least squares method in equation (\ref{eq:least_squares_law}).

\begin{table}
\label{page:table_inner_zeros}
\begin{center}
\caption{  $\#$ Zeros for Combined $C_{13}$ and $C_{23}$ Families}
\vspace{0.10in}
\begin{tabular}{| r  || r |r |}
\hline
{\bf Degree} & $\#$ {\bf Zeros} & {\bf Prediction} \\
\hline
5000 & 4 & 4.5 \\
10000& 6 & 6.5 \\
15000& 8 & 7.9 \\
20000& 9 & 9.1 \\
25000& 10 & 10.2 \\
30000& 11  & 11.2\\
35000& 12 & 12.1\\
40000& 13 & 12.9 \\
50000& 14 & 14.5 \\
60000& 17 & 15.8 \\
70000& 17 & 17.1 \\
\hline
\end{tabular}
\end{center}
\end{table}

\section{ Asymptotics Inside the Unit Disk}\label{section:asymptotics_inside}

The purpose of this section is to give a flavor of  the asymptotics of the partition
polynomials  needed to confirm their behavior already described.
Since the proofs are lengthy, we will be brief.

The generating function for the partition polynomials is a bivariate version
of the generating function for the partition numbers which gives
an integral version of the polynomials:
\[
P(x,u) = \sum_{n=1}^\infty F_n(x) u^n
=
\prod_{j=1}^\infty \frac{1}{1- xu^j},
\qquad
F_{n}(x)=
\frac{1}{2\pi i}\oint_{\left| u\right| =r}  \, \frac{P(x,u)}{u^{n+1}} \, du.
\]
We study $F_{n}(x)$ in each of the three regions ${\mathcal R}(k)$,
$k=1,2,3$ by using the Circle Method.
We begin by writing $F_n(x)$ as a contour integral over a circle
of radius $r$:
\begin{eqnarray*}
\lefteqn
{
F_n(x)
=
\frac{1}{2\pi i} \, \oint_{\left| u\right| =r}  \, \frac{P(x,u)}{u^{n+1}} \, du
} \\
&&=
e^{2 \pi n /M} \,
\sum  \int_{- \theta_{h,k}^\prime}^{ \theta_{h,k}^{\prime \prime}} \,
P(x, \exp[ 2 \pi i ( h/k + iz) ] ) \, e^{- 2\pi i n \, h/k} \, e^{- 2 \pi i n \phi} \, d \phi
\end{eqnarray*}
where the sum is over
rational numbers $h/k$, $(\theta_{h,k}', \theta_{h,k}'')$ are certain subarcs of
the circle that contain the point $e^{2 \pi i h/k}$. The choice of the radius $r$
and the order of the denominators $N$ of the fractions $h/k$ requires great care
-- this is the initial setup of the Circle Method using Farey fractions.
See  \cite{andrews}, \cite{ayoub}, or \cite{rademacher}.
Both the radius $r$ and the order of the Farey fractions $N$ depend on the
degree of the polynomial $n$ and the point $x$ inside the unit disk, see
\cite{boyer_goh_partition}.

%
%
%
%
%

The next step is to develop  the asymptotics of the integrals $I_{h,k}$
\[
I_{h,k}= 
 \int_{- \theta_{h,k}^\prime}^{ \theta_{h,k}^{\prime \prime}} \,
P(x, \exp[ 2 \pi i ( h/k + iz) ] )  \, e^{- 2 \pi i n \phi} \, d \phi.
\]
This requires a special approximation to the generating function
in a neighborhood of the rational points of the unit circle.
The proof of this is very lengthy and uses Dirichlet $L$-functions
and other techniques from analytic number theory. See our paper
\cite{boyer_goh_partition} for details.

For relatively prime integers $h<k\leq N$, we introduce the functions 
$Q_{h,k}(s)$  
that come from expanding $\ln[ P(x,e^{2 \pi i h/k + iz})]$ and are
given in the right plane $\Re(s) >1$
as the double series
\[
Q_{h,k}(s):=\sum_{m\geq 1}\sum_{l\geq 1}
\frac{x^{l}e^{2\pi ilmh/k}}{l}(lm)^{-s}.
\]
We establish that  
$Q_{h,k}(s)$ has an analytic continuation to $\mathbb C$
with a unique simple pole at $s=1$,
with residue at $s=1$ is
$ \frac{1}{k^{2}} \, \Li_2 (x^{k})$
and its evaluation at $0$ is
$Q_{h,k}(0)=
\frac{1}{2k}\ln (1-x^{k}) \,  + \, \sum_{l(k\nmid l)}\frac{x^{l}}{l}
\frac{1}{e^{-2\pi ilh/k}-1}.
$

\begin{theorem}\label{thm:factorization}
Let $h<k$ be relatively prime integers.
Fix  $|x|<1$, then there is a neighborhood of $u=e^{2\pi i h/k}$ inside
the unit disk where 
the generating function $P(x,u)$ has the factorization 
\[
P \left(x, e^{2 \pi i ( h/k + iz) }  \right) 
=
e^{ w_{h,k}} \, e^{ \Psi_k (z)} \, e^{ j_{h,k}(z)}, \quad \Re(z)>0,
\]
where
\begin{eqnarray*}
\Psi _{k}(z) &=& \frac{\Li_2 (x^{k})}{2\pi k^{2}} \, \frac{1}{z},\quad
j_{h,k}(z) = \frac{1}{2\pi i}\int_{-\frac{3}{4}+i\infty }^{-\frac{3}{4}+i\infty} \,
Q_{h,k}(s)  \, \Gamma (s)  \, (2\pi z)^{-s}   \, ds 
\\
w_{h,k} &=& 
 \frac{1}{2k}\ln (1-x^{k}) \,  + \, \sum_{l(k\nmid l)}\frac{x^{l}}{l}
\frac{1}{e^{-2\pi i  l  h/k}-1}
\end{eqnarray*}
\end{theorem}

Theorem \ref{thm:factorization} is crucial step that allows applying the Circle Method to the integrals $I_{h,k}$.
Let $f_{\rm max}(x) = \max[ f_1(x),f_2(x),f_3(x)]$
we state the results:
\begin{theorem}
(a) For every $k \geq 4$, there exists $\epsilon \in (0,1)$ so that
\[
|  I_{h,k} |
\leq
\frac{C}{k} \,
\frac{ f_{\max}(x)}{ \sqrt{\pi}} \, \frac{1}{n^{3/4}}
\,
\exp
\left[
\sqrt{n} \left(  (1- \epsilon) f_k(x) + f_{\max}(x) \right)
\right].
\]
(b)
For $k=1,2$, or $3$, on $\overline{ {\mathcal R}(k)} \cap U^+$:
\[
I_{h,k}
=
\frac{1}{  \sqrt{\pi}} \, \frac{1}{ n^{3/4}} \,
\left[ \frac{ \sqrt{ \Li_2(x^k)}  }  {k} \right]^{1/2} \,
\exp
\left(  2 \sqrt{n}  \frac{ \sqrt{ \Li_2(x^k)} }  {k} \right)
\,
\left( 1 + O \left( \frac{1} { \sqrt{n}} \right) \right)
\]
uniformly on compact subsets.
\end{theorem}

The above Theorem shows that
$I_{0,1}$ is the main contribution to $F_n(x)$ on region ${\mathcal R}(1)$,
$I_{1,2}$ the main contribution to $F_n(x)$ on region 
${\mathcal R}(2)$,
and
$I_{1,3}$ and $I_{2,3}$ are the main contributions
  to $F_n(x)$ on region ${\mathcal R}(3)$.
 In our last Theorem, we state this precisely.
 For convenience we introduce the notations
 $U^+$ for the upper open unit disk and $I_k$ for
\begin{equation}\label{eqn:integral_I_k}
I_k = 
\frac{1}{  \sqrt{\pi}} \, \frac{1}{ n^{3/4}} \,
\left[ \frac{ \sqrt{ \Li_2(x)}  }  {k} \right]^{1/2} \,
\exp
\left(  2 \sqrt{n}  \frac{ \sqrt{ \Li_2(x^k)} }  {k} \right).
\end{equation}

\begin{theorem}\label{theorem:polynomial_asymptotics}
(a) In $\overline{{\mathcal R}(1)} \cap U^+$:
$ 
F_n(x) = e^{w_{0,1}} I_1 + 
o_{K_1}
\left( o(I_1)  \right)
$
uniformly on compact subsets $K_1$.
\\
(b)
In $\overline{{\mathcal R}(2)} \cap U^+$:
$ 
F_n(x) = (-1)^n e^{w_{1,2}} I_2 + 
o_{K_2}\left( o(I_2)  \right)$
uniformly on compact subsets $K_2$.
\\
(c)
In $\overline{{\mathcal R}(3)} \cap U^+$:
$ 
F_n(x) =(e^{-2\pi in/3}e^{w_{1,3}}
+
e^{-2\pi  in2/3}e^{w_{2,3}}) \,   I_3 + 
o_{K_3}\left( I_3 \right)
$
uniformly on compact subsets $K_3$.
\end{theorem}

We can now describe how to normalize the polynomials $F_n(x)$ so
we can apply the Density Theorem (see Theorem \ref{theorem:density}).

We can apply Theorem \ref{theorem:polynomial_asymptotics}
to find the limits $\frac{1}{\sqrt{n}} \ln | F_n(x)|$. For example,
on region ${\mathcal R}(1)$,
\begin{eqnarray*}
\frac{1}{\sqrt{n}} \ln | F_n(x)|
&=&
\frac{1}{\sqrt{n}} \ln | I_1 ( e^{w_{01}} + o(1))|   \\
&=&
\frac{1}{\sqrt{n}} \ln\left| \frac{3}{4} \ln n - \ln(2 \pi) +
 \frac{1}{4} \ln[ \Li_2(x)]
+ 2 \sqrt{n} \sqrt{ \Li_2(x)}\right|
 + \frac{1}{\sqrt{n}} \ln| e^{w_{01}}+ o(1)|
\to
2 \Re \sqrt{ \Li_2(x)}.
\end{eqnarray*}

\section{Summary}

We gave an experimental exposition of our work on the zeros
of the partition polynomials that are a polynomial version
of partition numbers \cite{boyer_goh_partition}.
Intensive computational studies are required to understand their zeros
since the number of zeros inside the unit disk is proportional to the
square root of their degree. What is especially surprising is that three
families of such sparse zeros arise inside the disk. Furthermore, 
zeros for one of these families, the $C_{23}$ family,  do not occur until degree 13,000
and a second zero appears only at degree 60,000.

Our development of the asymptotics for these polynomials completely
determines the limiting behavior of the zeros both as points and their
density.  This gives an explanation of the rarity of zeros for $C_{23}$.

Among the many polynomial families given by a bivariate infinite
product generating function, we single out one related to plane
partitions.
It would be of great interest to find the zero asymptotics for the
  polynomial family $\{Q_n(x)\}$ where they are given
by the generating function:
\begin{eqnarray*}
\sum_{n=1}^\infty Q_n(x) u^n
=
\prod_{k=1}^\infty  \frac{1}{ (1-xu^k)^k}
\end{eqnarray*}
We write $Q_n(x) = \sum_{m=1}^n q_m(n) x^m$ as usual.
By \cite{andrews}, Chapter 11, or \cite{stanley_plane},
the coefficients
$q_m(n)$ count  the number of plane partitions of $n$
whose trace is $m$; that is, the sum of its diagonal entries
is $m$. Our last figure is a plot of the zeros of $Q_{100}(x)$.

\begin{figure}[h!]
\begin{center}
\includegraphics[height=5cm,width=5cm]{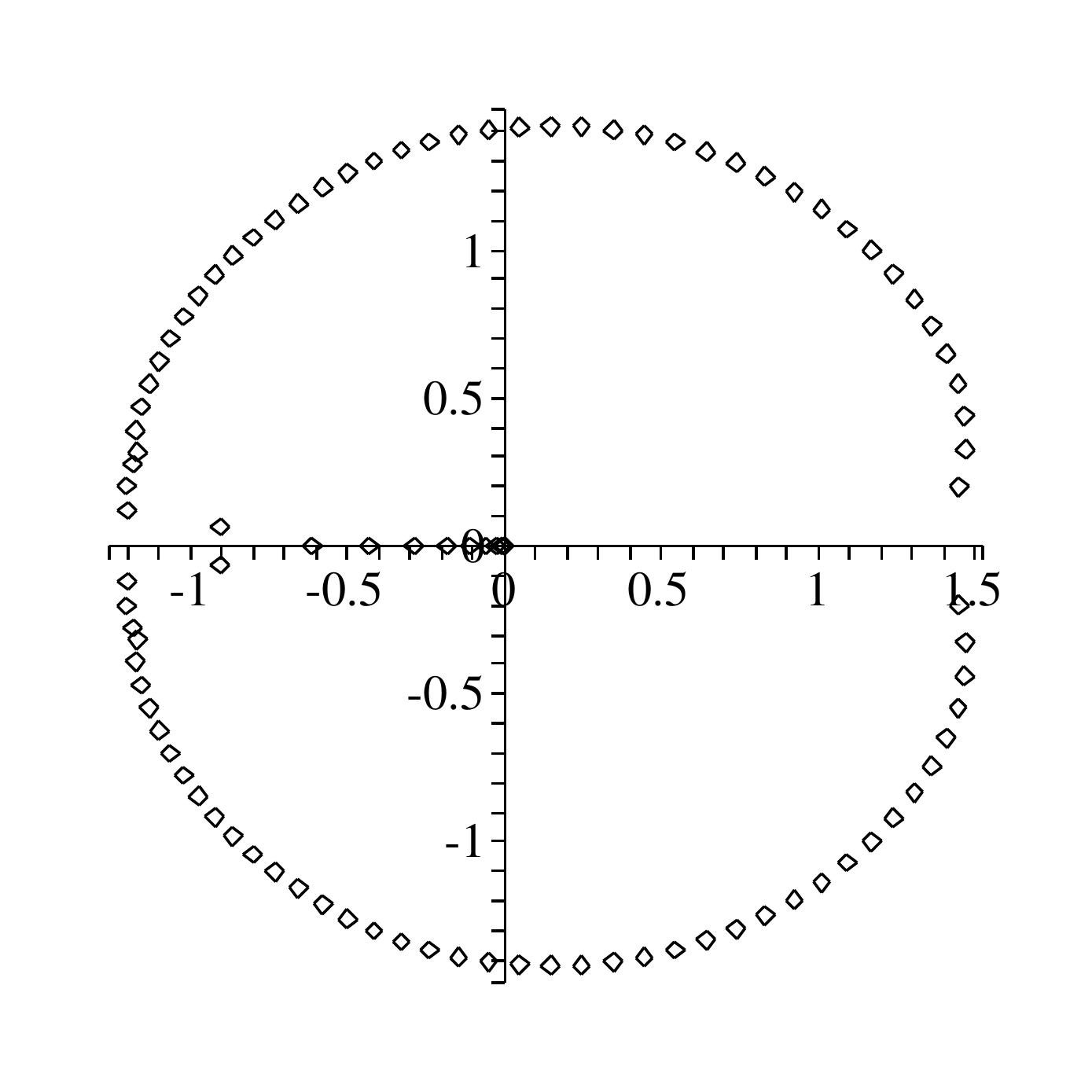} 
\caption{Plane Partition Polynomial Zeros of Degree 100}
\end{center}
\end{figure}

\vspace{0.3in}
\noindent
Department of Mathematics
\\
Drexel University
\\
Philadelphia, PA 19104
\\
{\sl email}: {\tt rboyer  at   math.drexel.edu}


\begin{thebibliography}{9}

\bibitem{andrews} George Andrews, 
{\sl The Theory of Partitions}, Addison-Wesley, 1976.


\bibitem{ayoub}  Raymond Ayoub, 
{\sl An Introduction to the Analytic Theory of
Numbers}, AMS Mathematical Surveys {\bf 10}, 1963.

\bibitem{bini} Dario Bini and G. Fiorentino,
Design, Analysis and Implementation of a Multiprecision Polynomial Rootfinder,
{\sl Numerical Algorithms} {\bf  23} (2000) 127-173.

\bibitem{euler}  Robert P. Boyer and William M. Y. Goh, On the zero attractor of Bernoulli
polynomials, {\sl  Advances in Applied Mathematics}
{\bf 38}
(2007) 97-132.

\bibitem{boyer_goh_partition}
Robert P. Boyer and William M. Y. Goh,
The Zero Attractor for Partition Polynomials, Manuscript.



\bibitem{hardy}
G. H. Hardy and S. Ramanujan, 
Asymptotic formulae in combinatorial analysis,
Proceedings of the London Math. Society {\bf } (1918)
75-115.


\bibitem{rademacher} 
Hans Rademacher, {\sl Topics in Analytic Number Theory},
Springer-Verlag, 1973.

\bibitem{saff_totik}
Edward B. Saff and Vilmos Totik,
{\sl
Logarithmic Potentials with External Fields},
Springer Verlag, 1997.

\bibitem{sokal} Alan Sokal,
Chromatic roots are dense in the whole complex plane,
Combinatorics, Probability \& Computing {\bf 13} (2004), 221-261.
 
 \bibitem{stanley_plane} Richard P. Stanley,
 The conjugate trace and trace of a plane partition,
 {\sl J. Combinatorial Theory} {\bf A14} (1973), 53-65.

\bibitem{varga}
Richard S. Varga and Amos J. Carpenter, 
Zeros of the partial sums of $\cos(z)$ and $\sin(z)$,
{\sl Numerical Algorithms} {\bf 25} (2000), 363-375. 


\end{thebibliography}
\end{document}